\documentclass{amsart}
\usepackage[utf8]{inputenc}
\usepackage{amssymb}
\usepackage{hyperref}
\usepackage{amsmath}
\usepackage{amsaddr}
\usepackage{appendix} 
\usepackage{graphicx}
\newcommand{\hR}{{^{\ast}\mathbb{R}}} 
\newcommand{\Ri}{{^{\ast}\mathbb{R}}_i} 
\newcommand{\Rf}{{^{\ast}\mathbb{R}}_f} 

\newcommand{\W}{\Ri[{\rm i}]}
\newtheorem{definition}{Definition}

\newtheorem{theorem}{Theorem}
\usepackage{url}

\title{An $\hR$-based perspective on solving ordinary differential equations}
\author{M.~Weber}
\address{Zuse Institute Berlin, Takustrasse 7, D-14195 Berlin, Germany}
\email{weber@zib.de}

\begin{document}

\begin{abstract}
 The real numbers, it is taught at universities, correspond to our idea of a continuum, although the hyperreal numbers are located ``in between'' the real numbers. The number $x + dx$, where $dx$ should be an infinitesimal number and $x$ real, is infinitesimally close to $x$ but ``infinitely'' far away from all other real numbers. Analogously: If $f'(x_0)$ and $f(x_0)$ are given for a differentiable function $f:\mathbb{R}\rightarrow\mathbb{R}$ at $x_0\in\mathbb{R}$, we can not determine $f(x)$ at {\em any} point $x\in \mathbb{R}$ different from $x_0$. These points seem to be ``infinitely'' far away. That is one conceptual problem of solving differential equations in numerical mathematics. In this article, we will present a numerical algorithm to solve very simple initial value problems. However, the change of paradigm is, that we will not ``leave'' the point $x_0$. Solving ordinary differential equations is like searching for ``recipes'' $f$.  Instead of trying to find these recipes for values $x\in\mathbb{R}$, we will learn them from special relations in the ``monad'' of $x_0$. MSC: 34A45, Keywords: non-standard analysis, infinitesimal roots, numerical treatment of ODEs 
\end{abstract}

\maketitle

\section{Introduction}\label{sec:Intro}


This article is a significantly (conceptually and mathematically) modified version of a preprint \cite{weber2020novel}.
Solving differential equations is usually understood as the searching for a differentiable function $f$ which solves a given equation. What does ``a function is  differentiable'' mean? In non-standard analysis the derivative of a function $f:\mathbb{R}\rightarrow \mathbb{R}$ is usually defined as \cite{AR60, BK16}:
\begin{equation}\label{eq:hrdiff}
\frac{df}{dx}(x)=st\Big(\frac{f(x+dx)-f(x)}{dx}\Big),
\end{equation}
where $dx$ is an infinitesimal number. The mapping $st(\cdot)$ is used to transform the fraction in (\ref{eq:hrdiff}), which is a hyperreal number, to a real number (the ``closest'' real number). 

The {\em standard part} $st(\cdot)$ is a piecewise constant function in terms of hyperreal numbers, which is not defined for infinite numbers. If the standard part of the fraction (\ref{eq:hrdiff}) exists and is independent from the choice of the infinitesimal number $dx$, then $f$ is said to be differentiable. In this context it is important that the function $f$ is at least continuous in $x$, which equivalently means that the nominator $f(x+dx)-f(x)$ is an infinitesimal number (otherwise, the fraction would be for sure infinite). Thus, for continuous functions the {\em fraction} is defined in the set of hyperreal numbers, but is not necessarily finite or independent from $dx$. This is similar to standard analysis, where ``being independent from the choice of the sequence $dx\rightarrow 0$''-condition restricts the set of functions which can be differentiated and can in turn be used as a possible solution of a differential equation. Many examples of real-valued, continuous functions $f:\mathbb{R}\rightarrow\mathbb{R}$ exist, which are not differentiable at any point $x\in \mathbb{R}$. One example - the Bolzano function - can be found in the script ``Functionenlehre'' of Bernhard Bolzano (arround 1830) and is available in \cite{Rus2004}. These are examples for functions where the fraction in (\ref{eq:hrdiff}) can be computed for every $dx$, but the standard part does not exist or is not independent from the choice of $dx$. 

\paragraph{\bf Alternative point of view.} $st(\cdot)$ marks the point in the differential equation where a transition from hyperreal numbers to real numbers occurs. In this article we move this transition point to the ``end'' of the calculation. Therefore, in our case, the functions $f$ can be seen as a mapping between hyperreal numbers. If we restrict ourselves to ``analytic'' functions, then a function value $f(x)$ is generated or approximated by a sequence of additions, subtractions and multiplications of hyperreal numbers including the input value $x$. Instead of looking for values $x$ and their function values $f(x)$, we will focus here on finding the {\em recipes} $f$ which solve an ordinary differential equation. The key point is that for such a recipe only operations are necessary that could also be carried out within an ideal of the ring of coefficients of a polynomial, so that we can restrict ourselves to the mapping between infinitesimal numbers. Such a mapping can be generalized to a mapping between fields. This article uses an algebraic analysis of hyperreal numbers. 

\paragraph{\bf Generalizations.} Other approaches try an analysis in terms of the even larger class of surreal numbers \cite{RuSw14}, where the  concept of ``limits'' of sequences is transferred to  surreal numbers. In this article we will go back to the 17th and 18{th} century mathematics of algebraic analysis \cite{PF21, G48} (before ``limits'' have been defined), but we will make use of hyperreal numbers \cite{AR60}. $\hR$ is the field of hyperreal numbers. The ring of infinitesimal numbers will be denoted as $\Ri\subset \hR$ and the ring of finite numbers will be denoted as $\Rf\subset \hR$. The fact that $\Ri$ is an ideal in $\Rf$ allows for algebraic analysis of polynomials. We will generalize the definition of derivatives to complex functions $f:\mathbb{C}\rightarrow \mathbb{C}$ by allowing for $dx \in\W$, where $\rm i$ is the imaginary unit. In this case, $st(\cdot)$ has to be applied to the real and complex part separately. The term ``infinitesimal number'' now also applies for elements of $\W$. The term ``finite numbers'' is used for $\Rf[{\rm i}]$. Every number in $\hR[{\rm i}]$ which is not finite is denoted as ``infinite''.  

Tab.~\ref{fig:scheme} gives a quick overview of our changes of the paradigm of how to solve a differential equation. The changes will be discussed in detail in this article.  
\begin{table}[h]
    \centering
    \begin{tabular}{|l|c|c|}
         \hline
         & well-known approach &  our novel approach\\
         \hline
         \hline
         search space & $f:\mathbb{C}\rightarrow \mathbb{C}$& $f\in\mathbb{C}[X], \beta\in\Ri[{\rm i}]$ \\
         \hline
         fixed & we search for one solution $f$ & we fix $\alpha$ in (\ref{eq:diff})\\
         \hline
         variation & $\forall x,\alpha$ in (\ref{eq:diff}) we want $F\in\Ri[{\rm i}]$ & $\forall f$ compute zeros $\beta$ of F\\
         \hline
         role of $f$ & global solution & numerical recipe at point $\beta$\\
         \hline
         philosophy & $st(F)=0$ & $st(x_0+\beta)=x_0$\\
         \hline
    \end{tabular}
    \caption{Conceptual changes to be discussed in this article}
    \label{fig:scheme}
\end{table}

\section{What is a solution?}\label{sec:definitions}
{\em In this section we show how we can move the transition point $st(\cdot)$ to the ``end'' of the calculation and what kind of functions and differential equations we want to discuss about. There are several related ways to define the term ``solution of a differential equation'' in this case.}

Solving a first-order ordinary differential equation is the same like searching for the zeros of a function $F(\frac{df}{dx},f,x)$. If there is a function $f$ such that $F(\frac{df}{dx},f,x)=0$ for all $x$, then $f$ is called {\em solution of the differential equation} \cite{WP20}. In this article, we will restrict ourselves to differential equations which are given by polynomials $F\in\mathbb{C}[Z,Y,X]$. Many differential equations can be transformed into such a polynomial\footnote{Rational functions $\widetilde{F}$ can also be transformed into $F=0$ with a polynomial $F$. Transcendent functions occurring in a differential equation can, furthermore, be approximated by rational functions in the domain of interest.}. So far, the concept for checking whether $f$ solves $F$ is like this: First, we have to apply the definition (\ref{eq:hrdiff}) for the computation of a complex-valued $\frac{df}{dx}$ and then we can check for $F=0$. Alternatively, we can also define the derivative of the function $f$ to be
\begin{equation}\label{eq:diff}
f'(x)=\frac{f(x+\alpha)-f(x)}{\alpha},
\end{equation}
where $\alpha\not= 0$ is an infinitesimal number $\alpha\in \W$. Now, we apply the rules of the function $st(\cdot)$ and exchange the application of $F$ and $st(\cdot)$. Thus, we search for functions $f$ with $F(st(f'),f,x)=st(F(f',f,x))=0$, i.e., in the well-known approach, we serach for functions $f$ with $F(f',f,x)\in\W$ for all $\alpha\in \W$ and $x\in\mathbb{C}$.

\paragraph{\bf Example 1.} Let us take the differential equation $F=Z-3X^2$. In standard and in non-standard analysis we would call $f(x)=x^3$ a solution of this equation, because $0 =\frac{df}{dx}(x)-3x^2$. In the given novel concept we first insert (\ref{eq:diff}) to the $Z$-component of the polynomial $F(Z,Y,X)$ and yield:
$$F(f'(x),f(x),x)=\frac{(x^3+3x^2\alpha+3x\alpha^2+\alpha^3)-x^3}{\alpha}-3x^2=\alpha\cdot(3x+\alpha).$$
The polynomial $f(x)=x^3$ is indeed a solution of the differential equation, but only in the following sense: For all $x\in\Rf[{\rm i}]$ and $\alpha\in \W$ the expression $F(f'(x),f(x),x)$ is an infinitesimal number. In terms of real numbers, $F$ can not be distinguished from being zero, but in terms of hyperreal numbers, $F$ is not zero (at least not for all $x\in\Rf[{\rm i}]$, in this example only at $x=-\frac{1}{3}\alpha$). The set $T_f(F)$ of zeros of $F$ for a given $f$ is very important in the next sections. Here is the new definition for the solution of a differential equation based on the discussed transformations:

\begin{definition}
According to (\ref{eq:diff}), a polynomial $f\in\mathbb{C}[X]$ is denoted as {\em hyper-solution of $F\in\mathbb{C}[Z,Y,X]$ in $D$}, if $F(f'(x),f(x),x)\in\W$ for all $x\in D$.
\end{definition}

The definition is even restricted to polynomials $f$, which has certain practical reasons.

\paragraph{\bf Is this point of view changing something?} This perspective does not open up any new possibilities {\em per se}. So far it is only said how we can determine for a function $f$ (found outside this framework) whether it is the solution of the differential equation $F$. Our answer would initially not differ very much from the answer of the standard analysis. Let us assume, however, that we have to determine $f$ numerically, since we do not know the algebraic solution of a differential equation, for example. In this case, such a solution is usually determined for selected points $x$. Recipes $f(x)$ are allowed that can be realized by a computer in a finite time (transcendent functions are approximated). The numerical perspective is about finite calculability and, thus, shares some concepts with finitism \cite{BM19}. Often:
\begin{itemize}
    \item[] For selected points $x \in \mathbb{C}$ recipes $f\in \mathbb{C}[X]$ are searched for, so that $f(x)$ approximates the solution of $F$ at $x$.  
\end{itemize}
Different recipes $f$ are selected for different points $x$ in a numerical routine. The actually determined solution $s:\mathbb{C}\rightarrow\mathbb{C}$ is then made up of these assignments of points $x$ and respective recipes $f$. The recipes are usually not specified explicitly, but result indirectly from a chaining of several numerical intermediate steps and possibly also by including several intermediate results on the ``path'' from a starting point $x_0$ to the point $x$. In addition, iterative solution methods can play a role, so that it would be laborious (but theoretically possible) to explicitly determine the currently used recipe $f$ for each point $x$. The new view described above allows this assignment to be reversed. So for given recipes $f$ that can be easily realized numerically now determine points $x$ at which the respective recipe represents a good approximation of the solution of $F$. And this is where the real trick comes in: 
\begin{itemize}
\item[] For selected recipes $f\in\mathbb{C}[X]$, points $\beta$ are searched for where equality $F(f'(\beta), f(\beta),\beta) = 0$ holds. 
\end{itemize}
The differential equation is therefore exactly fulfilled at these points. With this way of thinking, an error is immediately noticeable, which we see based on our standard analysis perspective: It is not possible to speak of a ``solution of a differential equation'' if $F=0$ is only realized at isolated points, but not ``in between''. This leads us back to the idea of the abstract. This article is a first step towards reversing the assignment between points $x$ and numerical recipes $f$ in a meaningful way. 

\paragraph{\bf Example 2.} Let us consider the equation $F=Z-Y$ with an initial value of $f(0)=1$, then in the set of polynomials $f\in\{\mathbb{C}[X], f(0)=1\}$ a solution of $F$ can not be found. The standard\footnote{This also would be a hyper-solution in $\mathbb{C}$, if we would allow for non-polynomials $f$.} solution would be given by $f(x)=\exp(x)$.  Let us enter the ``numerical way'' to solve an equation like $f'=f$. How do we calculate $\exp(x)$ in practice? The exponential $\exp(x)$ is approximated by a Taylor polynomial, e.g. $\exp(x)\approx f(x)=1+x+\frac{1}{2}x^2+\frac{1}{6}x^3$. If we insert the polynomial $f$ into the differential equation, the resulting $F$-polynomial is
\begin{equation*}
 F(f'(x),f(x),x) = -\frac{1}{6}x^3+\frac{1}{6}\alpha^2+\frac{1}{2}\alpha x+ \frac{1}{2}\alpha.
\end{equation*}
The leading monomial $-\frac{1}{6}x^3$ remains.  $F$ has the form $F=rx^n + \alpha\, G(x,\alpha)$ with $n>0$ and $r\not=0$. Surprisingly, this polynomial meets the requirements for a hyper-solution: $F$ is  infinitesimal for all $x\in \W=D$. The trick is to restrict $D$ to the set of infinitesimal numbers. This is the spirit of Taylor polynomials, they ``equal'' the functions and their derivatives only in ``one point'' - the center point. Any polynomial $f$ which leads to an $F$-polynomial with vanishing constant term, is a hyper-solution of $F$ in the set of infinitesimal numbers. For example, the same argumentation also holds for the polynomial $\bar{f}=1+x+\frac{1}{2}x^2+x^3$ which leads to an $F$-polynomial of the form: \begin{equation*}
\bar{F}=-x^3+\frac{5}{2}x^2+\alpha^2+3\alpha x+\frac{1}{2}\alpha.
\end{equation*}
$\bar{F}$ is infinitesimal for infinitesimal input values $x$. This polynomial has the form $\bar{F}=P(x)+\alpha\, G(x,\alpha)$, where $P$ does not have a constant term. 

\paragraph{\bf Example 3.} The above considerations show, that $F(f',f,x)$ in general has many different hyper-solutions $f(x)$, if we restrict $D$ to infinitesimal input values $x$, i.e., if we analyze $F$ in the ``monad'' \cite{G48} of the center point $0$. How can we find these hyper-solutions? Taking the example $F=(f')^2+f^2-1$ with an initial value $f(0)=y_0$ at the center point $0$. First, we fix the polynomial degree ($d_f=2$) and prepare the template polynomial $f(x)=ax^2+bx+c.$ Then we insert $f$ into the differential equation and apply $st(\cdot)$  to $F$ which provides  $$st(F)=a^2x^4+2abx^3+(4a^2+b^2+2ac)x^2+(4ab+2bc)x+b^2+c^2-1.$$ 
The initial value condition provides the constant term of $f$ which is $c=y_0$ (in our case we chose $y_0=0$). Now we proceed step by step through the monomials $x^k$ starting with $k=0$: The constant term of $st({F})$ should vanish, such that $f$ is a hyper-solution of $F$ for infinitesimal input values. Thus, the equation $b^2+c^2-1=0$ has to be solved. Together with $c=0$ this provides two possible solutions $b=1$ or $b=-1$. It is true, that the initial value problem leads to two different standard solutions ($\sin(x)$ and $-\sin(x)$). We will proceed with $b=1$. In principle we are already done: Any polynomial of the form $f(x)=ax^2+x$ leads to a vanishing constant term in $P$, where $F(f',f,x)=P(x)+\alpha\, G(x,\alpha)$. However, we can also take the higher monomials $x^k$ into account. For $k=1$, the equation $4ab+2bc=0$ with $b=1$ and $c=0$ has to be solved.\footnote{For the initial value condition $f(0)=1$, this equation would not add any further information to the search of the coefficients.} This equation provides $a=0$. Now all coefficients are fixed $a=0, b=1, c=0$. Interestingly, with these settings we get $st(F)=P=x^2$. Only one monomial remains, it turns out to be the leading one in $F$ also regarding the infinitesimal terms.  

\paragraph{\bf What kind of polynomials $F$ are possible?} Some (well-known) differential equations $F(Z,Y,X)$ are polynomials only in monomials $Z, X\! Z, Y,$ and in monomials of the form $X^n$ with a maximal degree $m$, see Tab.\ref{tab:examples1}. 
\begin{table}[h]
\centering
\begin{tabular}{|c|c|c|}
      \hline
     $F$ & $F(f',f,x)$ & one solution \cr
     \hline
     $XZ+Z-1$ & $(x+1)f'-1$ & $\ln(x+1)$\cr
     $Z-Y$ & $f'-f$ & $\exp(x)$ \cr
     $(X+1)Z+2Y$ & $(x+1)f'+2f$ & $(x+1)^{-2}$\cr
     $Z-3X^2-2X$ & $f'-3x^2-2x$ & $x^3+x^2$\cr
     \hline
\end{tabular}
\caption{\label{tab:examples1} Some examples of differential equations which allow for hyper-solutions $f$ such that $F=rx^n + \alpha\, G(x,\alpha)$. If $n$ has to be higher than the grade of $G$ with regard to $x$ (i.e. for $\alpha=0$), then $f$ coincides with Taylor polynomials of standard solutions.}
\end{table}

In these cases, the polynomial degree of $F$ with regard to $x$ (i.e. by setting $\alpha=0$) is the same as the polynomial degree $d_f$ of $f$ in $x$ (if this degree is at least $m$). We apply the above method for finding a hyper-solution, such that $F=rx^n + \alpha\, G(x,\alpha)$, where $n=d_f$ and a given initial condition $f(0)=y_0$ is valid. We get {\em one} condition for the coefficients from the initial value, $f(0)=y_0$. We get further $d_f$ conditions for the coefficients, because all terms $x^n$ with $n<d_f$ have to vanish in $F$. In total, we get $d_f+1$ conditions for $d_f+1$ coefficients. The resulting polynomials $f$ are Taylor approximations of the standard solution of $F$. 
\begin{table}[h]
\centering
\begin{tabular}{|c|c|c|}
      \hline
      $F$ &$F(f',f,x)$ & one solution \cr
     \hline
     $Z^2+Y^2-1$ & $(f')^2+f^2-1$& $\sin(x)$\cr
     $4(X+1)Z^2+Y^2-1$ & $4(x+1)(f')^2+f^2-1$& $\sin(\sqrt{x+1})$\cr
     $(X+1)^3Z-2Y$ & $(x+1)^3f'-2f$ &$\exp(-(x+1)^{-2})$\cr
     \hline
\end{tabular}
\caption{\label{tab:examples2} Some examples of differential equations which allow for hyper-solutions $f$ such that $F=P(x) + \alpha\, G(x,\alpha)$. }
\end{table}

In those cases, when $F(Z,Y,X)$ has terms like $Z^2$ or  like $X^3Z$, the polynomial degree of $F$ in $X$ is in general higher than the polynomial degree of $f$, see Tab.\ref{tab:examples2}.  If $f$ is a Taylor approximation of the standard solution, it can not be assured for such equations that $P(x)$ is of the form $rx^n$. However, it can be assured that the polynomial degree of $P$ is higher than the polynomial degree of $G$ with regard to $x$. These considerations lead to the following definitions:
\begin{definition}
A polynomial $f\in\mathbb{C}[X]$ is denoted as {\em hyper Taylor approximation of $F$}, if $F(f',f,x)$ has the form $F=rx^n + \alpha\cdot G(x,\alpha)$ with $G\in\mathbb{C}[X,Y]$, $n\in\mathbb{N}_{>0}$, and $r\in\mathbb{C}$, where the grade of $G$ with regard to $X$ is smaller than $n$.
\end{definition}
\begin{definition}
A polynomial $f\in\mathbb{C}[X]$ is denoted as {\em hyper local approximation of $F\in\mathbb{C}[Z,Y,X]$}, if $F(f',f,x)$ has the form $F=P(x) + \alpha\cdot G(x,\alpha)$ with $G\in\mathbb{C}[X,Y]$, where the grade of $G$ with regard to $X$ is smaller than the grade of $P$ and where $P$ does not have a constant term.
\end{definition}

As a summary of our findings we formulate the following
\begin{theorem}
Let $F\in\mathbb{C}[Z,Y,X]$ be a differential equation and $f\in\mathbb{C}[X]$ a complex polynomial, then the following statements hold.
\begin{itemize}
    \item[(i)] If and only if $F(f',f,x)=\alpha\cdot G(x,\alpha)$ with $G\in\mathbb{C}[X,Y]$, then f is a hyper-solution of $F$ in $\mathbb{C}$. 
    \item[(ii)] If $f$ is a hyper Taylor approximation and $r\not=0$, then $F(f',f,x)=0$ implies $x\in\W$.
    \item[(iii)] If $F(f',f,x)\in\mathbb{C}[x,\alpha]$ is lacking a constant term, then $f$ is a hyper-solution of $F$ in $\W$.
    \item[(iv)] Every hyper-solution of $F$ is also a hyper Taylor approximation. Every hyper Taylor approximation is also a hyper local approximation. 
\end{itemize}
\end{theorem}

Only the second statement needs to be shown. (Case 1) Assume, that $F=rx^n + \alpha\cdot G(x,\alpha)$ with $r\not=0$. Furthermore, assume that $x$ is a finite non-infinitesimal number.  Then $F$ is finite and non-infinitesimal. $x$ is not a zero of $F$. (Case 2) Assume, that $F=rx^n + \alpha\cdot G(x,\alpha)$ with $r\not=0$. Furthermore, assume that $x\not=0$ such that we can divide by $x^n$. The equation $rx^n + \alpha\cdot G(x,\alpha)=0$  is then equivalent to $r+\alpha\cdot \tilde{G}(1/x, \alpha)=0$, where $\tilde{G}\in\mathbb{C}[X,Y]$ is a suitable polynomial. This equation can not be solved by an infinite number $x$, because in this case $1/x$ would be infinitesimal and $st(r+\alpha\cdot \tilde{G}(1/x, \alpha))=r\not=0$. q.e.d.

\section{Transporting recipes}\label{sec:solving}
{\em In this section we will discuss how a found recipe $f$ which approximates the solution of a differential equation at some center point $x_0$ can be used to find a ``new recipe'' $f_{x_1}$ at a different center point $x_1$.}

The rules for hyper local approximations and for the hyper Talyor approximations in Sec. \ref{sec:definitions} are constructed in such a way, that $f$ is related to the Taylor series method for solving a differential equation at the center point $x_0=0$. Whenever one wants to shift this center point to a different value $x_0\not= 0$, then one has to analyze the differential equation $F(f'(x),f(x), x+x_0)$ instead. In this section we will see, that changing the center point has different algebraic effects on the hyper Taylor approximations. 

\paragraph{\bf Example 4.} Take the equation $F=xf'-1$. Trying to find the hyper Taylor approximation of this equation is like trying to expand $\ln(x)$ at $x_0=0$. If we would like to extend the logarithm at $x_0=1$, we have to solve the differential equation $F(f',f,x+1)$ which is equal to
\begin{equation}\label{eq:ln}
    \tilde{F}(f',f,x)=F(f',f,x+1)=(x+1)f'-1.
\end{equation}
Indeed, the polynomial $\tilde{f}(x)=x-\frac{1}{2}x^2$ leads to
$\tilde{F}(\tilde{f}',\tilde{f},x)=-x^2-\frac{1}{2}\alpha (x+1)$. Thus, $\tilde{f}(x)$ is a hyper Taylor approximation of $\tilde{F}$. This procedure leads to the following 

\begin{definition}
Let $F\in\mathbb{C}[Z,Y,X]$ define a differential equation $F(f',f,x)$, then $F_{x_0}(f',f,x):=F(f',f,x+x_0)$ is called the {\em differential equation $F$ at $x_0$}. 
\end{definition}

The last equation (\ref{eq:ln}) has shown, that $F_1=(x+1)f'-1$ is the differential equation $F=xf'-1$ at $x_0=1$. At this point, $f_1(x)=x-\frac{1}{2}x^2$ is a hyper Taylor approximation of this equation\footnote{Note, that indeed the polynomial $X-\frac{1}{2}X^2$ for $X=(x-x_0)$ is the quadratic Taylor expansion of $\ln(x)$ at $x_0=1$}.

Can we transport a hyper Taylor approximation to a new center point? The situation is trivial, if the polynomial $F(Z,Y,X)$ does not have an $X$-term, i.e., if transportation of $F$ to $F_{x_0}$ does not change the differential equation like for $F=f'-f$ or $F=(f')^2+f^2-1$. In these cases, any suitable $f$ for $F$ is also suitable for $F_{x_0}$. This leads to the following

\begin{definition}
Let $f_0\in\mathbb{C}[X]$ be a hyper Taylor approximation of $F_0$. If $f_{x_0}(x)=f_0(x)$ is a hyper Taylor approximation of $F_{x_0}$ for all values $x_0\in\mathbb{C}$, then the differential equation $F$ is denoted as {\em exponential-like solvable}. 
\end{definition}

We would expect that a transformation of the form $f_{x_0}(x)=f_0(x+x_0)$ is the ``correct'' transformation rule\footnote{In standard analysis: If a function $f(x)$ (like $\exp(x)$) solves $F$ at $0$, then $f(x+x_0)$ solves $F$ at $x_0$ (like $\exp(x+x_0)=\exp(x_0)\cdot\exp(x)$). In our setting this is not true anymore. Thus, our approach is really different from standard analysis.}. However, if a hyper Taylor approximation $f$ leads to the polynomial $F(f',f,x)=rx^n+\alpha G$, then a substitution of $x\rightarrow x+x_0$ in $F(f',f,x)$ turns a monomial $rx^n$ into a polynomial $r(x+x_0)^n$. A hyper Taylor approximation turns into a hyper local approximation. Only in the case of $r=0$ this transformation is vaild. This means, only for hyper-solutions this transformation can be applied. This is the case, whenever the polynomial $f$ solves $F$ in standard analysis. The polynomial $X^2$ which is a hyper Taylor approximation of $F=f'-2x$, can thus be ``transported'' to the hyper Taylor approximation $X^2+2x_0X+x_0^2$ at the point $x_0$. More precise: 

\begin{definition}
Let $f_0\in\mathbb{C}[X]$ be a hyper Taylor approximation of $F_0$. If $f_{x_0}(x)=f_0(x+x_0)$ is a hyper Taylor approximation of $F_{x_0}$ for all values $x_0\in\mathbb{C}$, then the differential equation $F$ is denoted as {\em polynomial-like solvable}. 
\end{definition}

What kind of transformation would ``keep'' the non-infinitesimal monomial $rx^n$ in $F$ of a hyper Taylor approximation? A transformation of the type $x\rightarrow \frac{x}{x_0}$ with $x_0\not= 0$ would do so. The equation $F=xf'-1$, is a corresponding example. It is not polynomial-like solvable, because a suitable polynomial $f$ for $F_0$ does not exist. It is also not exponential-like solvable. We will first introduce the definition and show its applicability afterwards. 

\begin{definition}
 Let  $f_1\in\mathbb{C}[X]$ be a hyper Taylor approximation of a differential equation $F_1\in\mathbb{C}[Z,Y,X]$ at $1$. If $f_{x_0}(x)=f_1(\frac{x}{x_0})$ is a hyper Taylor approximation of $F_{x_0}$ for all values $x_0\not=0\in\mathbb{C}$, then the differential equation $F$ is denoted as {\em logarithm-like solvable}. 
\end{definition}
     
\paragraph{\bf Proof.} It will be shown that $F=xf'-1$ is logarithm-like solvable. Thus, take a polynomial $f_1\in \mathbb{C}[X]$ which is a hyper Taylor approximation of $F_1=(x+1)f'-1$ (such a polynomial exists). Now it has to be shown, that $\tilde{f}=f_1(\frac{x}{x_0})$ is a hyper Taylor approximation of $F_{x_0}$. \begin{equation*}
\begin{split}
    F_{x_0}(\tilde{f}',\tilde{f},x) & = (x+x_0)\tilde{f}'-1 \\
    &= (x+x_0)\frac{\tilde{f}(x+\alpha)-\tilde{f}(x)}{\alpha}-1 \\
    &= (x+x_0)\frac{f_1(\frac{x+\alpha}{x_0})-f_1(\frac{x}{x_0})}{\alpha}-1\\
    &= (x+x_0)\frac{f_1(\frac{x}{x_0}+\frac{\alpha}{x_0})-f_1(\frac{x}{x_0})}{\alpha}-1\\
    &= (x+x_0)\frac{1}{x_0}\frac{f_1(\frac{x}{x_0}+\frac{\alpha}{x_0})-f_1(\frac{x}{x_0})}{\frac{\alpha}{x_0}}-1\\    
    &= \big(\frac{x+x_0}{x_0}\big) \big(f_1'\big(\frac{x}{x_0}\big) +R_1\big)-1\\
    &= \big(\frac{x}{x_0}+1\big) f_1'\big(\frac{x}{x_0}\big)-1 + R_2\\
    &= (z+1)f_1'(z)-1 +R_2.
    \end{split}
\end{equation*}
Changing the infinitesimal quantity $\alpha$ to $\alpha/x_0$ changes the value of the derivative in (\ref{eq:diff}). However, it only changes the value up to an infinitesimal difference $R_1$ (the grades in $x$ are not changed). Since $R_2=R_1(x+x_0)x_0^{-1}$ and $(x+x_0)x_0^{-1}$ is finite, $R_2$ is also infinitesimal. Note, that $(z+1)f_1'(z)-1+R_2$ has the correct form for a hyper Taylor approximation, because $f_1$ is a hyper Taylor approximation of $F_1$, $R_2$ is infinitesimal (in polynomial form), and the leading monomial is transformed from $rx^n$ into $\frac{r}{x_0^n}x^n$ by $x\rightarrow z$.  q.e.d.

\paragraph{\bf Example 5.} With a similar calculation one can show that
\begin{equation}\label{eq:recip}
    F(f',f,x)=xf'+n\cdot f,
\end{equation}
with a natural number $n$, is logarithm-like solvable\footnote{For $n=1$ a hyper Taylor approximation of $F_1$ is given by $f(x)=1-x+x^2$.}.  The standard analysis solution of this equation would be $x^{-n}$. 

\paragraph{\bf Example 6.} For the complicated example $\exp(-1/x^2)$ of standard analysis, where the Taylor series at $x_0=0$ does not coincide with the function itself\footnote{after continuation with $\exp(-1/0^2):=0$}, the differential equation is
\begin{equation}\label{eq:taylor}
    F(f',f,x)=2f - x^3f'.
\end{equation}
The equation (\ref{eq:taylor}) does not look like polynomial-like solvable. However, one would have to check this for all hyper Taylor approximations. The polynomial $f=0$ is a hyper-solution of this differential equation, i.e., a hyper Taylor approximation. The corresponding (transported) polynomial is $f_{x_0}(x)=f(x+x_0)=0$ which also is a hyper Taylor approximation of the transported differential equation. If this $f=0$ is the only hyper Taylor approximation of (\ref{eq:taylor}), then it is polynomial-like solvable.  \\

The three definitions of Sec.~\ref{sec:solving} provide possible recipes to ``solve'' ordinary differential equations. If we can find approximates at every center point $x_0$, i.e., suitable polynomials $f_{x_0}$ for every $F_{x_0}$, then we can (at least locally) ``solve'' ordinary differential equations, because we know how the set of ``numerical recipes'' looks like locally for every number $x_0\not=0\in\mathbb{C}$. How can we turn hyper-solutions in $\W$ into standard solutions of differential equations in $\mathbb{C}$?

\section{Numerically motivated choice of the recipe} \label{sec:numerics}
{\em After reading this section it will be clear that the main problem of a numerical treatment of ordinary differential equations is given by the fact, that local approximations at a center point $x_0$ of a solution of $F$ are in general not valid solutions for any center point different from $x_0$.}

This section will follow the usual concept of numerical mathematics (like ``walking along the real axis''). We construct solutions of the differential equations with standard numerical tools \cite{DB13} like step size control and like adjustment of the polynomial degree. In all of the following cases, $s(t)$ is not a hyper Taylor approximation of $F$, unless there exists a hyper-solution of $F$ in $\mathbb{C}$. It is just a good numerical approximation of a standard solution of $F$. We will present the changes of paradigm in Sec.~\ref{sec:algebra}.

Although the procedure in Sec.~\ref{sec:solving} provides hyper Taylor approximations of a given differential equation $F$ for all points $x_0\in\mathbb{C}$, it does not seem to be satisfactory in terms of ``solving the differential equation''. We would expect that there is only {\em one function for all points $x_0\in\mathbb{C}$} instead of a set of functions (approximates) for each point. How to glue these local approximations together to yield an approximate global solution? In numerical mathematics, differential equations are treated in terms of {\em initial value problems}. In addition to $F(f',f,x)=0$, we further define an initial condition $f(x_0)=y_0$ to be satisfied. 
Solving an initial value problem like this in the context of this article, would mean to restrict the set of polynomials $f$ to a certain subset $f\in \mathbb{P}\subset \mathbb{C}[X]$, which meets the initial value condition $\mathbb{P}=\{f\in \mathbb{C}[X]; f(x_0)=y_0\}$. As an example look at $F=f'-f$. If we have found a polynomial which is a hyper Taylor approximation $f$ of this equation, then every multiple $\lambda\cdot f$ also is a hyper Taylor approximation in this special case. Only if we additionally ask for $f(0)=1$, then solely polynomials with constant part $1$ are valid. One example is $f(x)=1+x+\frac{1}{2}x^2\in\mathbb{P}$.  

The transportation mechanisms described in Sec.~\ref{sec:solving} not only have to transport the center point, but they also have to transport the initial value condition $f(x_0)=y_0$ to the new center point $t$. 

\begin{definition}
If an initial value problem is given by a feasible set $\mathbb{P}$ of polynomials and by a differential equation $F$, then we call $F_t$ together with the feasible set $\mathbb{P}_t=\{f\in\mathbb{C}[X]; f(X-t)\in\mathbb{P}\}$ the {\em initial value problem at $t\in\mathbb{C}$}.
\end{definition}

\underline{First idea} (locally solve the transported problems): For every $t\in\mathbb{C}$ we ``approximately solve'' the initial value problem at $t$. Let $f_t(x)$ denote the ``approximate solution'' of the initial value problem at $t$. Then, one would expect, that the function $s(t)=f_t(0)$ ``approximately solves'' the corresponding initial value problem. Note, that $f_t(x)$ has the center point $t$ and meets the required (transported) initial value condition $f(x_0)=y_0$. Let us check this na\"ive way: 

\paragraph{\bf Good Example.} This example shows how this transported solution of initial value problems works:
\begin{itemize}
    \item Take the example of the initial value problem $F=f'-2x$, with $f(0)=0$ and the hyper-solution $f(x)=x^2$.
    \item We have to find a feasible hyper-solution for $F_1$ and $\mathbb{P}_1$. The equation is $F_1=f'-2x-2$ and we need a  polynomial $f_1$ with $f_1(-1)=0$. The polynomial $f_1(x)=(x+1)^2=x^2+2x+1$ is feasible, because $f_1\in\mathbb{P}_1$, and $f_1$ is a hyper-solution of $F_1$.
    \item For the equation $F_2=f'-2x-4$, we need a polynomial $f_2$ with $f_2(-2)=1$. A corresponding hyper-solution is $f_2(x)=(x+2)^2=x^2+4x+4$.
    \item Thus, $f(x)=x^2$ is a hyper-solution at $t=0$, $f_1(x)=x^2+2x+1$ is a hyper-solution at $t=1$, and $f_2(x)=x^2+4x+4$ is a hyper-solution at $t=2$. 
\end{itemize} The three polynomials satisfy the (transported) initial value condition $f(0)=1$. In this situation, we would expect, that a ``solution'' $s(t)$ of the initial value problem is given by $s(0)=f(0)=0$, and $s(1)=f_1(0)=1$, and $s(2)=f_2(0)=4$, which coincides with $s(t)=t^2$. 

\paragraph{\bf Bad Example.} The next example shows, that this transportation mechanism is not valid for hyper Taylor approximations in general:
\begin{itemize}
    \item Take the example of the initial value problem $F=f'-f$, with $f(0)=1$. A feasible hyper Taylor approximation is $f(x)=1+x+\frac{1}{2}x^2$.
    \item We have to find a feasible hyper Taylor approximation for $F_1$ and $\mathbb{P}_1$. The equation is the same $F_1=f'-f$, but in this situation we search for a polynomial $f_1$ with $f_1(-1)=1$. The polynomial $f_1(x)=2+2x+x^2$ is feasible, because $f_1\in\mathbb{P}_1$, and $f_1$ is a hyper Taylor approximation of $F_1=F$. Note, that $f_1$  is just a multiple of $f$.
    \item For the equation $F_2=f'-f$, we need a polynomial $f_2$ with $f_2(-2)=1$. A corresponding hyper Taylor approximation is $f_2(x)=1+x+\frac{1}{2}x^2$.
    \item  Thus, $f(x)=1+x+\frac{1}{2}x^2$ is a hyper Taylor approximation at $t=0$, $f_1(x)=2+2x+x^2$ is a hyper Taylor approximation at $t=1$, and $f_2(x)=1+x+\frac{1}{2}x^2$ is a hyper Taylor approximation at $t=2$.
\end{itemize}
The three polynomials satisfy the (transported) initial value condition $f(0)=1$. In this situation, we would expect, that a ``solution'' $s(t)$ of the initial value problem is given by $s(0)=f(0)=1$, and $s(1)=f_1(0)=2$, and $s(2)=f_2(0)=1$, which does not coincide with $\exp(t)$. This is a bad approximation, because $f_t$ only hyper-solves the differential equation locally and the initial value condition (at a different position $x_0\not=t$) is out of this infinitesimal range. For higher order polynomials, like $f=\sum_{n=1}^{10}\frac{x^n}{n!}$, this procedure provides better estimates\footnote{For $f=\sum_{n=1}^{10}\frac{x^n}{n!}$ the estimates are $s(1)=\frac{45360}{16687}$ and $s(2)=\frac{2835}{383}$ which is close to the corresponding values of $\exp(t)$.} of $\exp(t)$. Asking  for polynomials $f$
with ``infinte'' grade to solve the initial value problem (including problems of convergence of Taylor series) is not the spirit of this article. \\

\underline{Second idea} (controlling the grade of the polynomials): For a numerical treatment of initial value problems we could further restrict the set of polynomials to ``better fits'', if we do not want to deal with ``infinite grades''. A hyper Taylor approximation $f_t\in\mathbb{P}_t$ of an initial value problem $\{F, \mathbb{P}\}$ at $t$ is a good approximation, if e.g. 
\begin{equation}\label{eq:goodapp}
|st(F_t(f',f,x))|<\epsilon, \text{ for all } |x+\frac{1}{2}t|\leq \frac{1}{2}|t|.
\end{equation}
Such a condition assures, that $f_t$ is a good numerical approximation of the differential equation in the whole ``interval'' $[-t,0]$ and that $f_t$ has the correct initial value. This additional condition (\ref{eq:goodapp}) further restricts the set\footnote{Also an empty set can be the result of this restriction.} of possible polynomials $f_t$. It's like a discretization-based adjustment of the grade of the polynomial. Again $s(t)=f_t(0)$ is the resulting numerical solution.\\

\underline{Third idea} (step size control): An alternative approach using a step size control on $t$ is  possible, too. Let us start with a differential equation $F$ and an initial condition $\mathbb{P}$. A hyper Taylor approximation $f$ of this initial value problem is also a good numerical solution of the differential equation in a certain ``region'', such that there is e.g. a $\Delta >0$ with $|st(F(f',f,x))|\leq \epsilon$ for all $|x|\leq \Delta$. We can evaluate $f$  at any point $t\in\mathbb{C}$ with $|t|\leq \Delta$ for a good approximation of the ``solution''. Select one value $t$ and $f(t)=y_0$. Then we proceed with the same argumentation for the next step  by replacing the differential equation $F$ with $F_t$ and the initial value condition $\mathbb{P}$ with the new condition  $f(0)=y_0$. In each step $k$ of this procedure we get a step size $t_k$ and a hyper Taylor approximation $f^{(k)}$. The numerical solution is given by $s(\sum_{k=1}^n{t_k})=f^{(n)}(0)$.

\paragraph{\bf Numerical Example.} As an example take the initial value problem $F=(x+1)\cdot f'-1$ with $f(0)=0=s(0)$. 
\begin{itemize}
    \item[$(1)$] A hyper Taylor approximation is given by $f^{(1)}(x)=x-\frac{1}{2}x^2$. We select a step size: $t_1=\frac{1}{2}$. This leads to $f^{(1)}(\frac{1}{2})=\frac{3}{8}$. 
    \item[$(2)$] The next initial value problem to be solved is $F=(x+\frac{3}{2})\cdot f'-1$ with $f^{(2)}(0)=\frac{3}{8}$.  A hyper Taylor approximation is given by $f^{(2)}=\frac{3}{8}+\frac{2}{3}x-\frac{2}{9}x^2$ which has been found by the logarithm-like solution procedure. Again selecting $t_2=\frac{1}{2}$. This yields $f^{(2)}(\frac{1}{2})\approx 0.65278$.
\end{itemize} 
The result is $s(0)=0$, $s(0.5)=0.375$, and $s(1)\approx 0.65278$. The ``true solution'' would be: $\ln(1+0)=0, \ln(1+0.5)\approx 0.40547, \ln(1+1)\approx 0.69315$.

\begin{figure}[htb]
    \centering
    \includegraphics[width=0.4\textwidth]{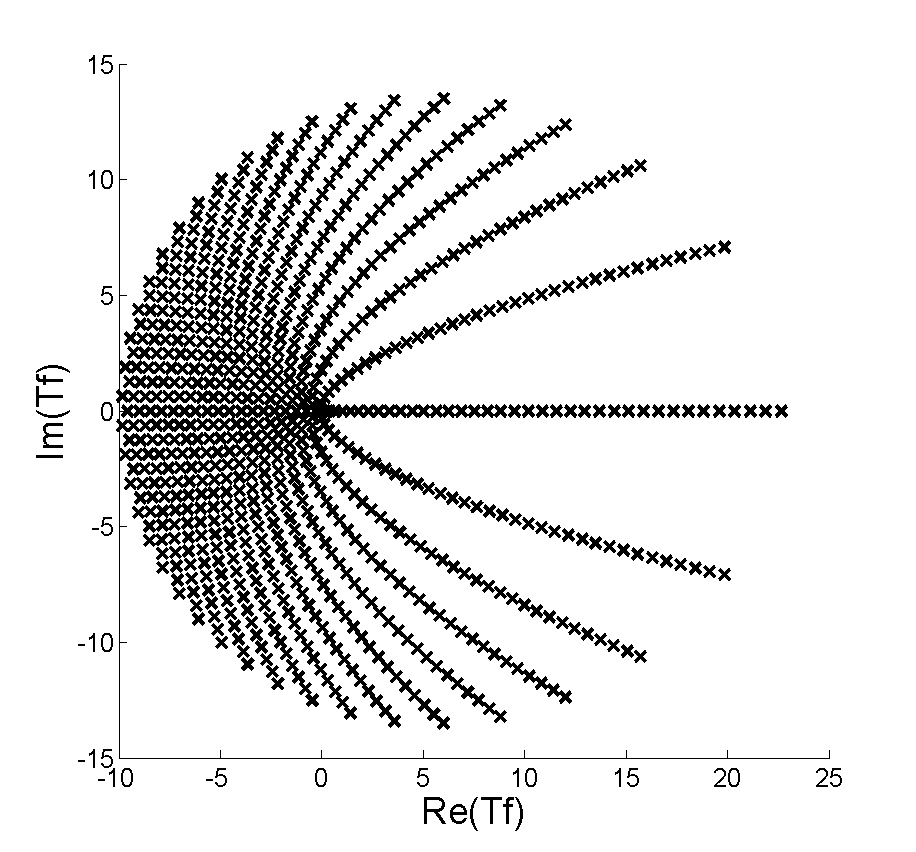}\\
    \includegraphics[width=0.4\textwidth]{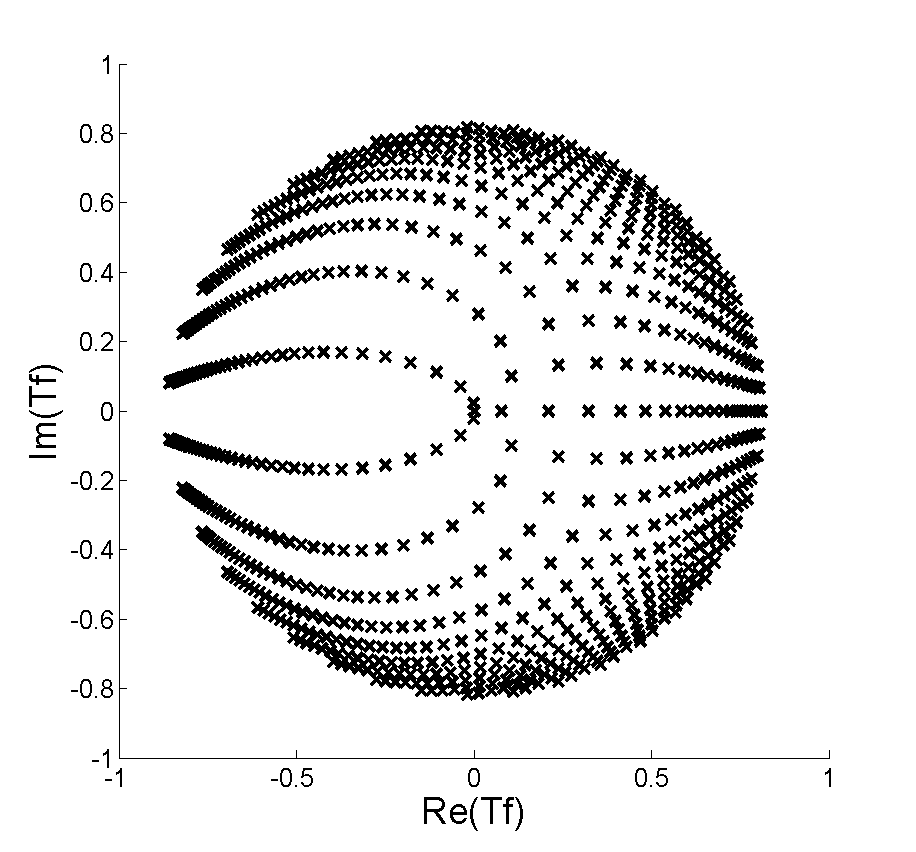}
    \includegraphics[width=0.4\textwidth]{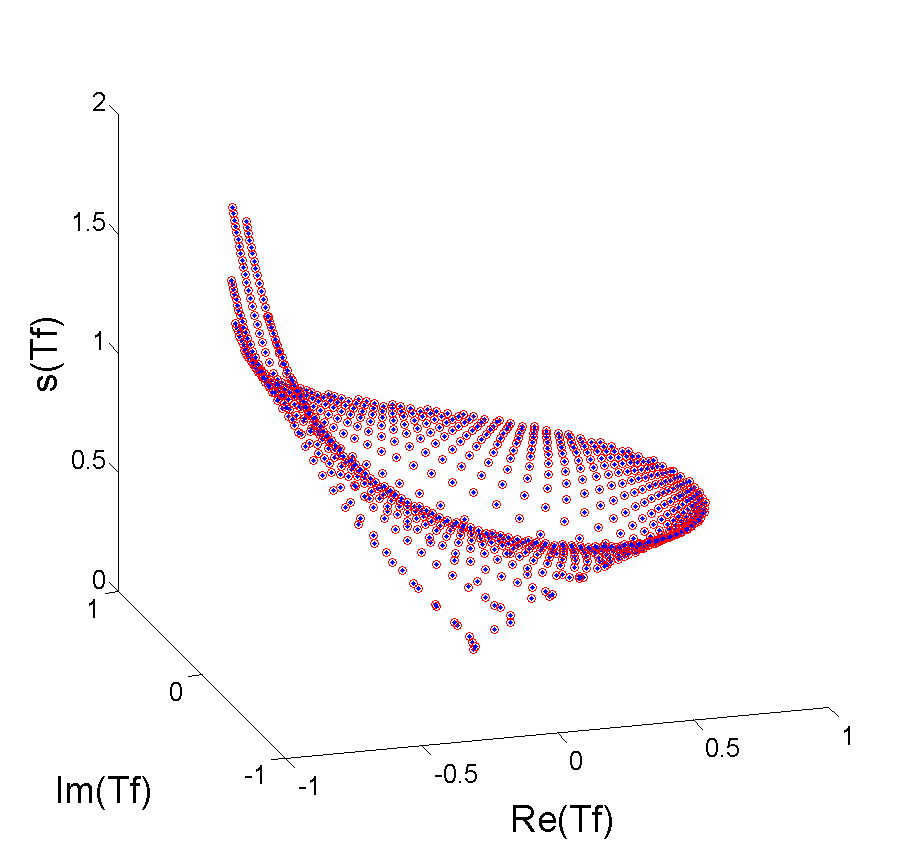}
    \caption{These figures show some results that have been found by {\sf insola} with $\alpha^\ast=0.001$ and a maximal polynomial degree of $n=40$. Top: The roots $T^\ast\subset \mathbb{C}$ for the initial value problem $F=f'-f$ and $f(0)=1$. Bottom Left: The roots $T^\ast\subset\mathbb{C}$ for the initial value problem $F=(x+1)f'-1$ and $f(0)=0$. Bottom right: Again, the same roots like in the left figure, but for every point $t^\ast\in \mathbb{C}$ the absolute value of $s(t^\ast)$ is shown on the $z$-axis.} 
    \label{fig:roots}
\end{figure}

\section{Algebraically motivated choice of the recipe} \label{sec:algebra}
{\em In this section we really want to find ``solutions'' that lead to $F=0$. The ``price'' we have to pay is, that for given functions $f$ this equation is only realizable for a finite set of infinitesimal roots $\beta$. All these values $\beta$ are indistinguishable from $0$ with regard to real numbers and that's the key!}

The problem of numerical methods described in the last section is the following: We want to get away from the center point $x_0=0$. Once we leave this point, the selected polynomials only approximate the solution of the equation. Thus, we should stay at the center point or in its ``monad''.

Hyper local approximations can not be distinguished from a solution of the differential equation in this case.  Through the glasses of real numbers, moving away by infinitesimal steps is not different from looking at $x_0=0$. If we insert one hyper-solution $f$ in $D=\W$, then the corresponding polynomial $F(f',f,x)$ only has a finite set of infinitesimal roots with  $F(f'(\beta), f(\beta),\beta)=0$ denoted as $T_f(F)$. A root $\beta \in T_f(F)$ provides a class of sequences, which exactly solves the equation for the given polynomial $f$ and the given class of sequences $\alpha$ for computing $f'$.  This means, that the solution $s$ should be constructed in the following way: Given a hyper-solution $f$ of $F$  in $D=\W$, we compute the infinitesimal roots of $F(f',f,\cdot)$ leading to the set $T_f(F)$. For every $\beta\in T_f(F)$ we define $s(\beta)=f(\beta)$, because it is $f(\beta)$ which is inserted into $F$ to provide a zero of $F$. We repeat this procedure for every hyper-solution $f$. Here it has to be said that only a relation $s:\W\leftrightarrow \W$ (not necessarily a function) is constructed in this way, because the sets of roots can intersect for different polynomials $f$. The interesting point is, that by computing the roots of $F$ for {\em every} given hyper-solution $f$, we determine at which points $\beta$ the algebraic recipe $f$ is a valid solution method. This procedure is very different from trying to find {\em one} solution function $f$ for which $F$ is infinitesimal on $\mathbb{C}$. Now we have recipes in the ideal $\W$. How can we make these recipes visible in $\mathbb{C}$?   
\begin{figure}
    \centering
    \includegraphics[width=0.4\textwidth]{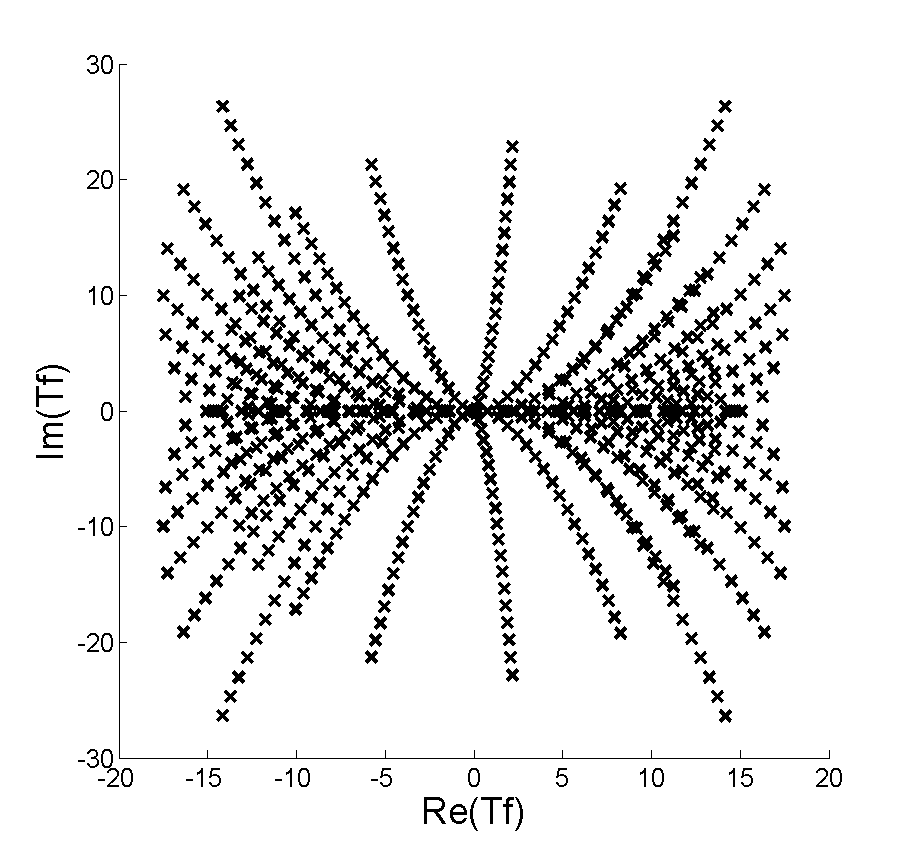}
    \includegraphics[width=0.4\textwidth]{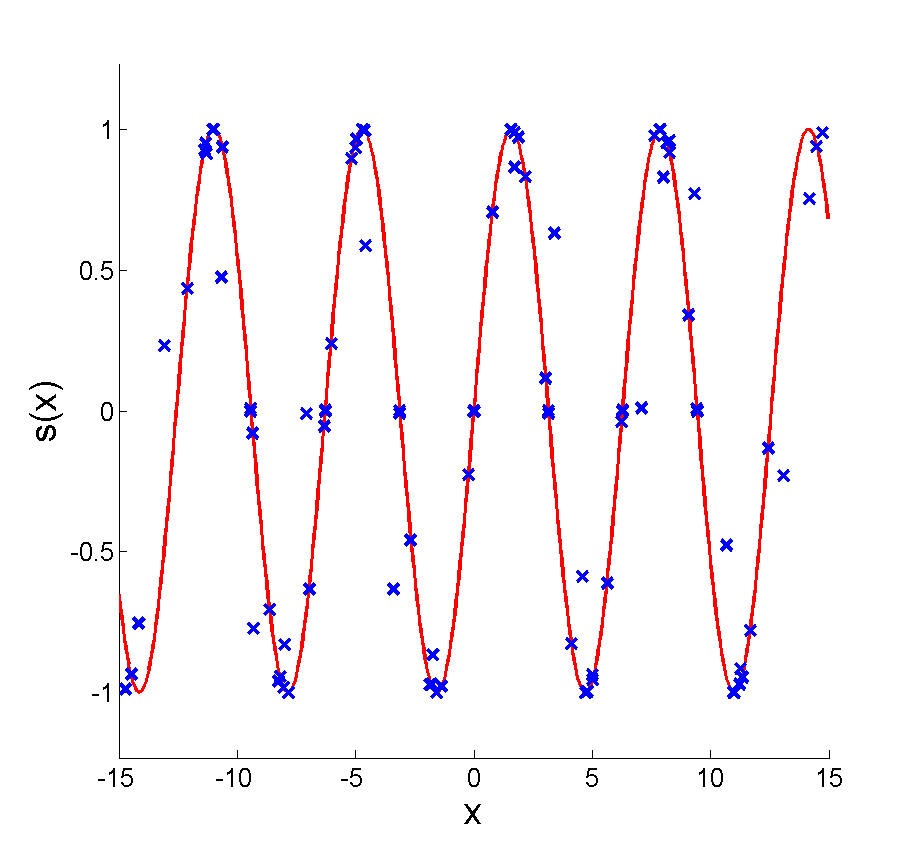}
    \includegraphics[width=0.4\textwidth]{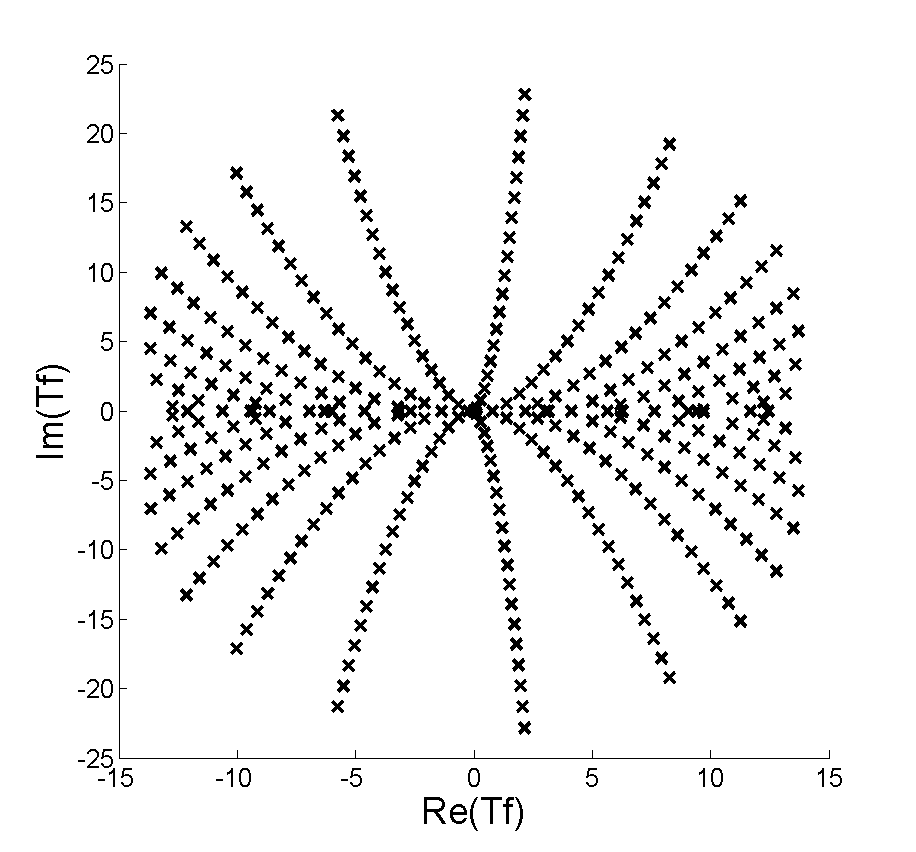}
    \includegraphics[width=0.4\textwidth]{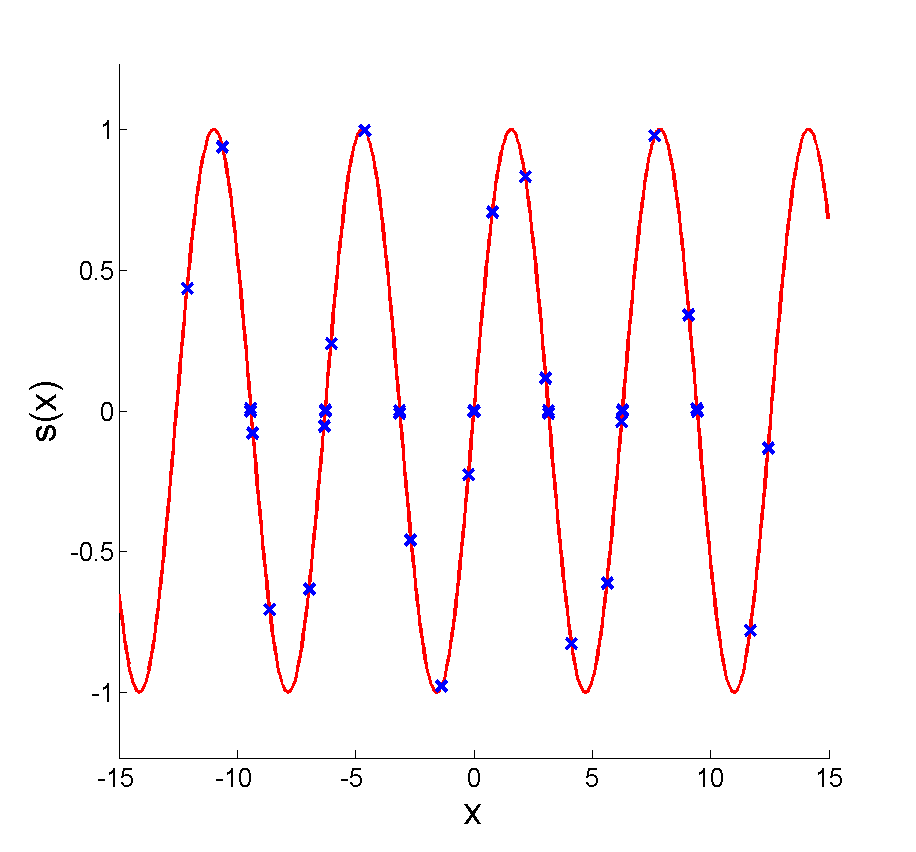}
    \caption{These figures show results that have been found by {\sf insola} with $\alpha^\ast=0.001$ and a maximal polynomial degree of $n=40$ with hyper local approximations instead of hyper Taylor approximations. Top Left: The roots $T^\ast\subset \mathbb{C}$ for the initial value problem $F=(f')^2+f^2-1$ and $f(0)=0$. Top right: Restriction to the real valued roots (on the x-axis) and the expected real valued results of $\sin(x)$ (red curve). Bottom Left: After eliminating the non-infinitesimal roots. Bottom right: Comparing the real-valued results from {\sf insola} with $\sin(x)$ after eliminating the non-infinitesimal roots.} 
    \label{fig:roots2}
\end{figure}

\paragraph{\bf Numerical experiments.} The described relation $s$ is numerically not accessible, because we apply infinitesimal numbers, which are not represented in numerical routines. However, if we want to find an {\em approximate} representation of this relation with non-infinitesimal numbers, then we maybe simply replace $\alpha$ with a very small real number $\alpha^\ast$ in the above considerations. Here comes the {\em infinitesimal solution algorithm} ({\sf insola}). Repeat for all grades $n$:
\begin{itemize}
    \item[1.] We first compute a hyper Taylor approximation $f$ (with $r\not=0$) of the initial value problem for a fixed polynomial grade $n$. 
    \item[2.] Then we determine $F(f'(x),f(x),x)\in\Rf[{\rm i}][x]$ using (\ref{eq:diff}). $F$ has roots only in the set of infinitesimal numbers. 
    \item[3.] In the expression for $F$ we replace $\alpha$ with a small real value $\alpha\rightarrow \alpha^\ast$ and yield $F^\ast\in\mathbb{C}[x]$. 
    \item[4.] Then, we compute all roots of $F^\ast$, which are finite complex numbers $t^\ast\in T_f^\ast$. 
    \item[5.] For every number $t^\ast\in T_f^\ast$ we plot the relation $(t^\ast, f(t^\ast))\in s$.
\end{itemize}    
A MATLAB$^{\sf TM}$-code that can be used to visualize and do experiments with the different initial value problems is in the Appendix \cite{MATLAB:2017b}. In this code we replaced the search for hyper Taylor approximations by computing the Taylor polynomials (of the known solutions) directly. Without knowing the solution one could apply the method described in Example 3 in Sec.\ref{sec:definitions} in order to construct these polynomials. The roots of $F$ found by the algorithm for two different initial value problems are shown in Fig.~\ref{fig:roots} (on the top $\exp(x)$ and in the bottom row $\ln(x+1)$). We yield the same picture like on the bottom left (however with different function values for the points), if we solve the differential equation $F=xf'+2f$ which is also logarithm-like solvable according to (\ref{eq:recip}).  

In order to compare these representations with the expected function values in $\mathbb{C}$, the {\sf insola} points are plotted as blue points, whereas, the results from standard analysis are plotted as red circles (bottom right, Fig.~\ref{fig:roots}). The results from {\sf insola} coincide with our expectations about the standard solutions of the corresponding differential equations. 
\begin{figure}
    \centering
    \includegraphics[width=0.4\textwidth]{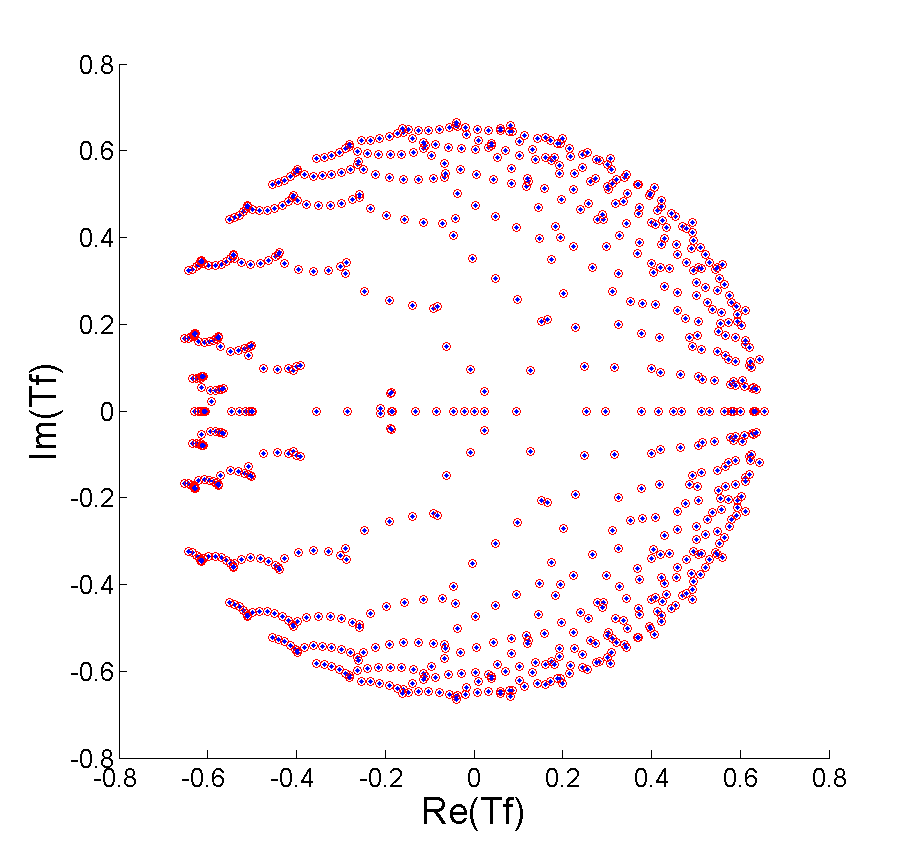}
    \includegraphics[width=0.4\textwidth]{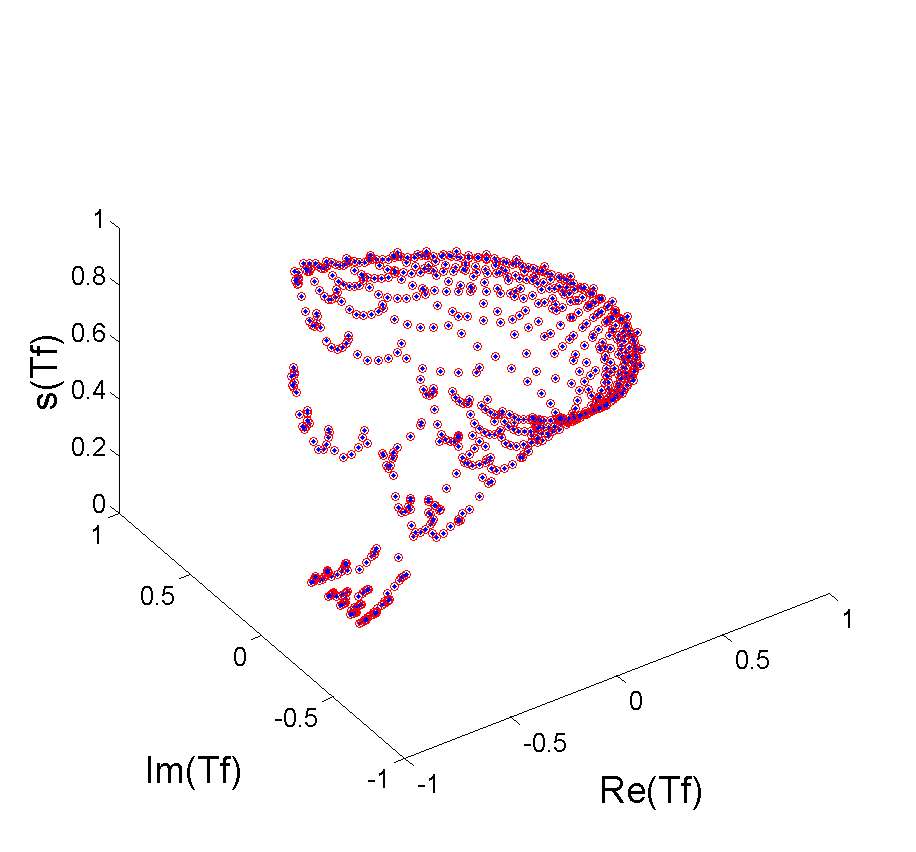}
    \caption{The {\sf insola} algorithm also finds the expected solution $\exp(-\frac{1}{(x+1)^2})$ of $F=2f-(x+1)^3f'$ after elimination of the non-infinitesimal roots. Left: The complex valued roots. Right: The roots versus the absolute value of the proposed solution $s$ on the z-axis (blue dots). The values coincide with the expected function values (red circles).} 
    \label{fig:roots3}
\end{figure}

\begin{figure}
    \centering
    \includegraphics[width=0.4\textwidth]{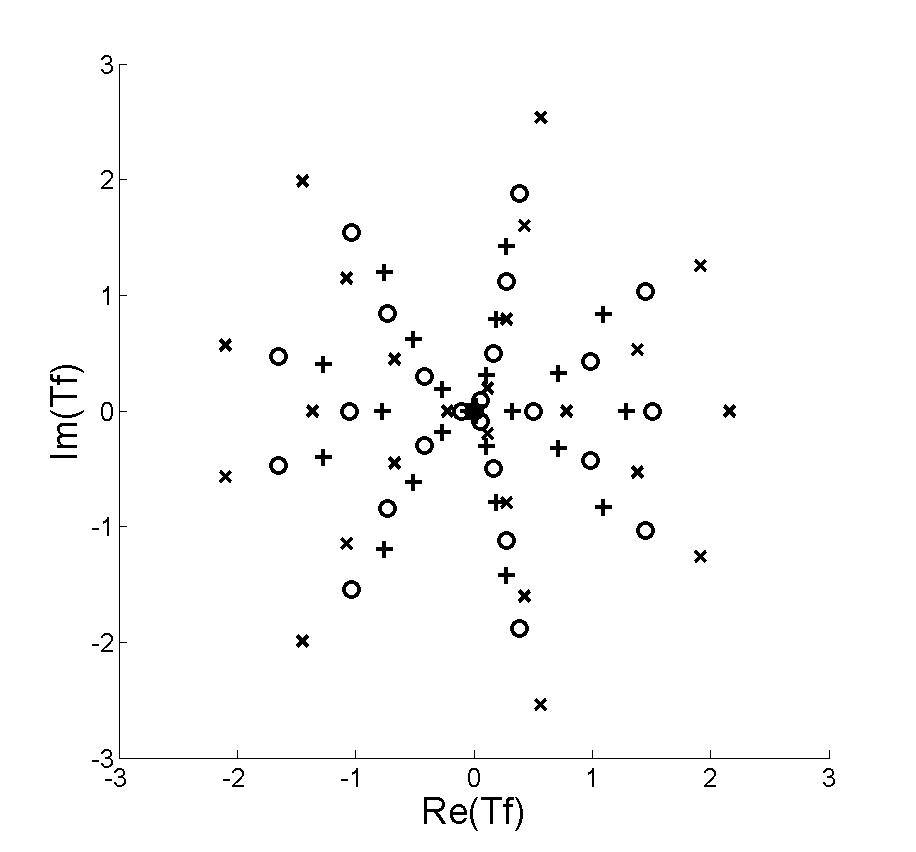}
    \includegraphics[width=0.4\textwidth]{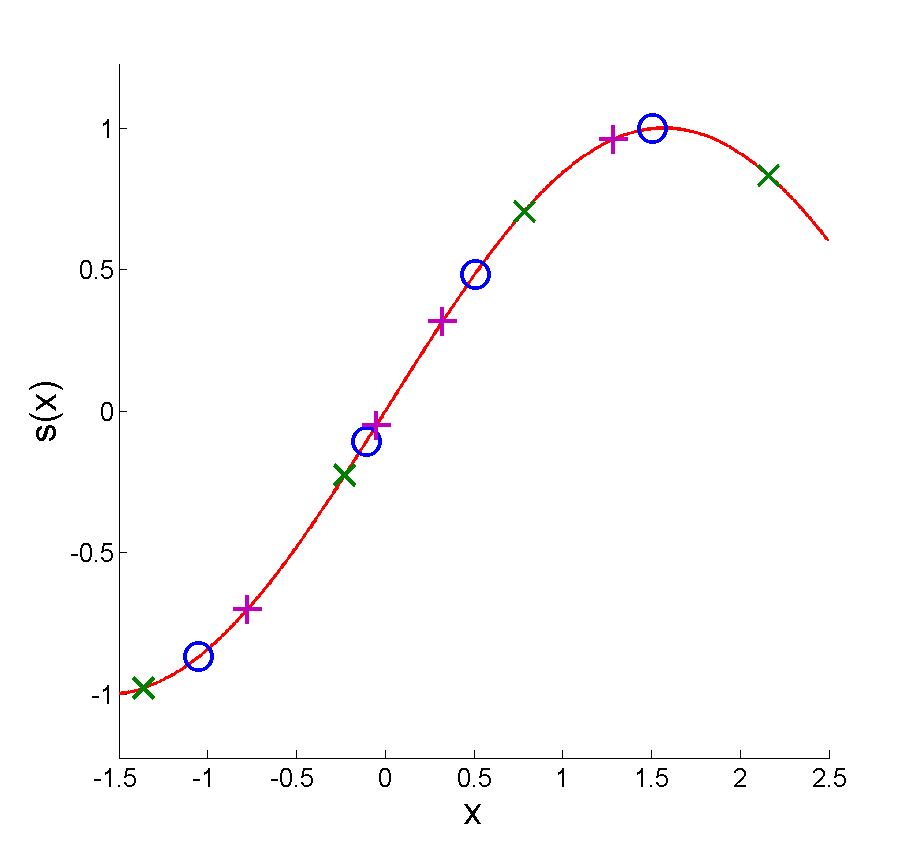}
    \caption{These plots show different representations (i.e. different values for $\alpha^\ast$, the symbol x means $\alpha^\ast=0.001$, circle means $0.0001$, and crosses means $0.00001$) for the approximation of the $\sin(x)$-function. Left: Changes of the complex roots. The smaller $\alpha^\ast$, the more the roots converge to zero. Right: Changes of the real-valued roots and comparison with the expected sin-function values.} 
    \label{fig:different}
\end{figure}

\paragraph{\bf Hyper local approximations.} In order to test the algorithm {\sf insola} for the case of hyper local approximations, we let it run for $F=(f')^2+f^2-1$ and the initial condition $f(0)=0$ and the additional condition that the coefficient of $x$ in $f$ is positive. We expect the solution $\sin(x)$. In the top right plot in Fig.~\ref{fig:roots2}, the real valued roots $x\in T_f^\ast$ are plotted versus the actual value of the approximated solution $s(x)=f(x)$ with blue crosses. However, these crosses do not coincide with the red curve, which shows the expected values ($x$ versus $\sin(x)$). Here another reason shows up, why {\sf insola} is based on hyper Taylor approximations. Hyper local approximations can lead to finite non-infinitesimal roots of $F$. The polynomial $f$, however, is only valid hyper-solution within the range of infinitesimal numbers. This means, the finite non-infinitesimal values are out of the region where $F$ is hyper-solved by $f$. If we insert a finite real value $\alpha^\ast$ into $F$ in step (3.) of {\sf insola} for numerical reasons, then the set $T_f^\ast$ which has to be constructed in the fourth step  will always consist of finite complex numbers. We can not distinguish between numbers that stem from ``turning $\alpha$ into a real number'' or from the finite non-infinitesimal roots of $F$. If we could  distinguish these two cases, then we could sort out the non-infinitesimal roots. The algorithm that we propose to sort out non-infinitesimal roots from $T^\ast_f(F)$ uses the assumption that $F$ is of the form $F(x,\alpha)=P(x)+\alpha G(x,\alpha)$, with a polynomial $P$. The non-infinitesimal roots of $T_f(F)$ are assumed to be close to the roots of $P$, which can be accessed by  setting $\alpha\rightarrow 0$ in $F(x,\alpha)$. We applied this method to the example of Fig.~\ref{fig:roots2} and indeed end up with roots and approximated values, which coincide with $\sin(x)$ (in the real valued roots and -not shown- also in the complex roots).  It also works for other examples, see Fig.\ref{fig:roots3}. The strategy is always the same: In the set of infinitesimal numbers, the hyper local approximations are the substitute for a ``solution'' of $F$. Now, we only regard the infinitesimal roots of $F$. By the replacement of $\alpha$ with $\alpha^\ast$ we can represent or approximate the scaled-up roots of $F$ in the complex plane. In this way, this procedure turns out to provide an approximate solution also in the set of complex numbers. Note that this ``scaling-up'' trick also works for different choices of $\alpha^\ast$, see Fig.~\ref{fig:different}

\section{Conclusion}
Usually {\em solving a differential equation} means to search for a function $f:\mathbb{C}\rightarrow\mathbb{C}$ such that $F(f'(x),f(x), x)$ is an infinitesimal number for all $x\in \mathbb{C}$. If we reduce (change) the input domain of $f$ to infinitesimal numbers, then {\em finding $f$} is not a problem anymore. There are many possible polynomials $f$ which are hyper-solutions of $F$ in this domain. In this new approach, the choice of $f$ defines the set $T_f(F)$ of points at which $F=0$. By this, we define our solution $s$. The interesting observation is, that it is possible to scale-up this relation $s$ into the set of complex numbers providing good approximations to the standard solutions of $F$. Maybe this new paradigm (``Search for $x$-dependent polynomials to compute $s$ in $\W$'') allows for a different perspective onto the existence and uniqueness (is $s$ a function?) of the solution of differential equations, in general.
 
\appendix

\section*{Appendix: Matlab Code}
\begin{verbatim}

syms x    
syms a   
syms f
syms fder %derivative of f

% standard solution for comparison with s() in the plots

fs=log(x+1);    %only leading monomial in P
%fs=exp(x);     %only leading monomial in P
%fs=sin(x);     %polynomial P
%fs=sin(sqrt(x+1));    %polynomial P
%fs=(x+1)^(-2); %only leading monomial in P
%fs=x^3+x^2;    %polynomial P=0 
%fs=exp(-(x+1)^(-2));  %polynomial P

% define small alpha and maximal grade of the polynomials 

alpha=0.001;
grades=40;


figure(1);
hold on;
figure(2);
hold on;

for i=2:grades
    % Taylor polynomials
    
    f=taylor(fs,i);

    % differential equation (select the correct one)
    
    dgl=(x+1)*fder-1;      %for log(x+1)
    %dgl=fder-f;           %for exp(x) 
    %dgl=fder^2+f^2-1;     %for sin(x)
    %dgl=4*(x+1)*fder^2+f^2-1; %for sin(sqrt(x+1)) 
    %dgl=(x+1)*fder+2*f;   %for (x+1)^(-2)
    %dgl=fder-3*x^2-2*x;   %for x^3+x^2
    %dgl=(x+1)^3*fder-2*f; %for e^(-(x+1)^(-2))
    
    %prepare F(f',f,x)
    F=expand(subs(dgl, fder, (subs(f,x, x+a)-f)/a));
    
    %numerical approximation: insert a finite alpha
    Fn=expand(subs(F,a,alpha));
    
    %roots of the polynomial Fn
    C=sym2poly(Fn);
    Tf=roots(C);
    Tf=unique(Tf);
    
    %insert a=0 to identify non-infinitesimal roots
    F0=expand(subs(F,a,0));
    
    %roots of the polynomial F0 (to exclude from Tf)
    C0=sym2poly(F0);
    Tf0=roots(C0);
    Tf0=unique(Tf0);
    
    %excluding non-infinitesimal roots
    for j=1:length(Tf0)
        if(Tf0(j)~=0)
          [val, ind]=sort(abs(Tf-Tf0(j)));
          Tf=Tf(ind(2:end));
        end
    end
    
    for j=1:length(Tf)
      % for the roots that are real plot s as a "graph"
      if(isreal(Tf(j)))
        figure(1);
        plot(Tf(j), subs(f,x,Tf(j)),'.b'); 
        plot(Tf(j), subs(fs,x,Tf(j)),'ro');
      end
      
      % for complex valued roots r plot Re(r), Im(r) against the 
      % absolute value of s(r)
      figure(2);
      plot3(real(Tf(j)),imag(Tf(j)), abs(subs(f,x,Tf(j))),'.b');
      plot3(real(Tf(j)),imag(Tf(j)), abs(subs(fs,x,Tf(j))),'ro');
    end
end


\end{verbatim}

\bibliographystyle{amsplain}
\bibliography{mycitations}
\nocite{*}
\end{document}